\documentclass[10pt]{amsart}
\usepackage{amsfonts}
\usepackage{ifthen}
\usepackage{amsthm}
\usepackage{amsmath}
\usepackage{graphicx}
\usepackage{amscd,amssymb,amsthm}
\usepackage{graphicx}
\usepackage{epstopdf}
\usepackage{hyperref}

\newcounter{minutes}
\divide\time by 60
\newcounter{hours}
\multiply\time by 60 \addtocounter{minutes}{-\time}

\newtheorem{theorem}{Theorem}

\newcommand{\real}{\operatorname{Re}}

\keywords{Starlike functions; radius of starlikeness; Mittag-Leffler expansions; $q$-Bessel functions; zeros of $q$-Bessel functions; Laguerre-P\'olya class of entire functions.}
\subjclass[2010]{30C45, 30C15, 33C10}

\title{Bounds for radii of starlikeness of some $q$-Bessel functions}

\author[\.{I}. Akta\c{s}]{\.{I}brah\.{I}m Akta\c{s}}
\address{Department of Mathematical Engineering, Faculty of Engineering and Natural Sciences, G\"{u}m\"{u}\c{s}hane University, G\"{u}m\"{u}\c{s}hane, Turkey}
\email{aktasibrahim38@gmail.com}
\author[\'A. Baricz]{\'Arp\'ad Baricz$^{\bigstar}$}
\address{Department of Economics, Babe\c{s}-Bolyai University, Cluj-Napoca, Romania}
\address{Institute of Applied Mathematics, \'Obuda University, Budapest, Hungary}
\email{bariczocsi@yahoo.com}

\thanks{$^{\bigstar}$The research of \'A. Baricz was supported by a research grant of the Babe\c{s}-Bolyai University for young researchers with project number GTC-31777.}

\begin{document}

%==================================================================
\def\thefootnote{}
\footnotetext{ \texttt{File:~\jobname .tex,
          printed: \number\year-0\number\month-\number\day,
          \thehours.\ifnum\theminutes<10{0}\fi\theminutes}
} \makeatletter\def\thefootnote{\@arabic\c@footnote}\makeatother
%==================================================================

\maketitle

\begin{abstract}
In this paper the radii of starlikeness of the Jackson and Hahn-Exton $q$-Bessel functions are considered and for each of them three different normalization are applied. By applying Euler-Rayleigh inequalities for the first positive zeros of these functions tight lower and upper bounds for the radii of starlikeness of these functions are obtained. The Laguerre-P\'olya class of real entire functions plays an important role in this study. In particular, we obtain some new bounds for the first positive zero of the derivative of the classical Bessel function of the first kind.
\end{abstract}

\section{Introduction}

Let $\mathbb{D}_r$ be the open disk $\{z\in\mathbb{C}:\left|z\right|<r\}$ with radius $r>0$. Let $\mathcal{A}$ denote the class of analytic functions $f:\mathbb{D}_r\rightarrow\mathbb{C},$ which satisfy the normalization conditions $f(0)=f^{\prime}(0)-1=0$. By $\mathcal{S}$ we mean the class of functions belonging to $\mathcal{A}$ which are univalent in $\mathbb{D}_r$ and let $\mathcal{S}^*$ be the subclass of $\mathcal{S}$ consisting of functions which are starlike with respect to origin in $\mathbb{D}_r$. The analytic characterization of this class of functions is $$\mathcal{S}^*=\left\{f\in\mathcal{S}:\real\left(\frac{zf^{\prime}(z)}{f(z)}\right)>0\text{ for all }z\in\mathbb{D}_r\right\}.$$ The real number $$r^{\ast }(f)=\sup \left\{ r>0 \left|\real\left(\frac{zf^{\prime }(z)}{f(z)}\right) >0 \;\text{for all }z\in\mathbb{D}_r\right.\right\}$$ is called the radius of starlikeness of the function $f$. Note that $r^{\ast }(f)$ is the largest radius such that the image region $f(\mathbb{D}_{r^*(f)})$ is a starlike domain with respect to the origin. For more information about starlike functions we refer to Duren's book \cite{Duren} and to the references therein.

Now, consider the Jackson and Hahn-Exton $q$-Bessel functions which are defined as follow:
$$J_{\nu}^{(2)}(z;q)=\frac{(q^{\nu+1};q)_{\infty}}{(q;q)_{\infty}}\sum_{n\geq0}\frac{(-1)^{n}\left(\frac{z}{2}\right)^{2n+\nu}}{(q;q)_{n}(q^{\nu+1};q)_{n}}q^{n(n+\nu)}$$ and $$J_{\nu}^{(3)}(z;q)=\frac{(q^{\nu+1};q)_{\infty}}{(q;q)_{\infty}}\sum_{n\geq0}\frac{(-1)^{n}z^{2n+\nu}}{(q;q)_{n}(q^{\nu+1};q)_{n}}q^{\frac{1}{2}{n(n+1)}},$$ where $z\in\mathbb{C},\nu>-1,q\in(0,1)$ and $$(a;q)_0=1,	(a;q)_n=\prod_{k=1}^{n}\left(1-aq^{k-1}\right),	 (a,q)_{\infty}=\prod_{k\geq1}\left(1-aq^{k-1}\right).$$
It is known that the Jackson and Hahn-Exton $q$-Bessel functions are $q$-extensions of the classical Bessel function of the first kind $J_{\nu}$. Clearly, for fixed $z$ we have $J_{\nu}^{(2)}\left((1-z)q;q\right)\rightarrow{J}_{\nu}(z)$ and $J_{\nu}^{(3)}\left((1-z)q;q\right)\rightarrow{J}_{\nu}(2z)$ as $q\nearrow1.$ The readers can find comprehensive information on the Bessel function of the first kind in Watson's treatise \cite{Wat} and properties of Jackson and Hahn-Exton $q$-Bessel functions can be found in \cite{ismail1,ismail2,koelink,Koornwinder} and in the references therein. The geometric properties of some special functions (like Bessel, Struve and Lommel functions of the first kind) and their zeros in connection with these geometric properties were intensively studied by many authors (see \cite{aktas2,aktas1,mathematica,publ,lecture,BCDT,bdoy,BKP,bsk,bos,samy,basz,basz2,BY,brown,ismail3,todd,wilf}). Also, the radii of starlikeness and convexity of some $q$-Bessel functions were investigated in \cite{BDM}. In some of the above mentioned papers it was shown that the radii of univalence, starlikeness and convexity are actually solutions of some transcendental equations. In addition, it was shown that the obtained radii satisfy some inequalities. On the other hand, it was proved that the radii of univalence and starlikeness of some normalized Bessel and Struve functions of the first kind coincide. Most of above papers benefited from some properties of the positive zeros of Bessel, Struve and Lommel functions of the first kind. Also, the Laguerre-P\'olya class $\mathcal{LP}$ of real entire functions, which consist of uniform limits of real polynomials whose zeros are all real, was used intensively (for more details on the Laguerre-P\'olya class of entire functions we refer to \cite{BDM} and to the references therein). Motivated by the earlier works, in this study our aim is to obtain some lower and upper bounds for the radii of starlikeness of some normalized $q$-Bessel functions. The results presented in this paper complement the results of \cite{BDM} about the radii of starlikeness and extend the known results from \cite{aktas1} on classical Bessel functions of the first kind to $q$-Bessel functions. As in \cite{BDM} we consider three different normalized forms of Jackson and Hahn-Exton $q$-Bessel functions which are analytic in the unit disk of the complex plane. Because the functions $J_{\nu}^{(2)}(.;q)$ and $J_{\nu}^{(3)}(.;q)$ do not belong to $\mathcal{A}$, first we consider the following six normalized forms as in \cite{BDM}. For $\nu>-1$,
$$f_{\nu}^{(2)}(z;q)=\left(2^{\nu}c_{\nu}(q)J_{\nu}^{(2)}(z;q)\right)^{\frac{1}{\nu}}, \nu\neq0$$
$$g_{\nu}^{(2)}(z;q)=2^{\nu}c_{\nu}(q)z^{1-\nu}J_{\nu}^{(2)}(z;q),$$
$$h_{\nu}^{(2)}(z;q)=2^{\nu}c_{\nu}(q)z^{1-\frac{\nu}{2}}J_{\nu}^{(2)}(\sqrt{z};q),$$
$$f_{\nu}^{(3)}(z;q)=\left(c_{\nu}(q)J_{\nu}^{(3)}(z;q)\right)^{\frac{1}{\nu}}, \nu\neq0$$
$$g_{\nu}^{(3)}(z;q)=c_{\nu}(q)z^{1-\nu}J_{\nu}^{(3)}(z;q),$$
$$h_{\nu}^{(3)}(z;q)=c_{\nu}(q)z^{1-\frac{\nu}{2}}J_{\nu}^{(3)}(\sqrt{z};q),$$ where $c_{\nu}(q)=(q;q)_{\infty}\big/(q^{\nu+1};q)_{\infty}$. In this way, all of the above functions belong to the class $\mathcal{A}$. Of course there exist infinitely many other normalization for both Jackson and Hahn-Exton $q$-Bessel functions, the main motivation to consider these six functions is the fact that their limiting cases for Bessel functions appear in literature, see for example \cite{brown} and the references therein.

\section{Bounds for the radii of starlikeness of  some normalized $q$-Bessel functions}

In this section our aim is to present some tight lower and upper bounds for the radii of starlikeness of the above mentioned six normalized forms of the Jackson and Hahn-Exton $q$-Bessel functions. In particular, we obtain some known and new bounds for the first positive zero of the first derivative of the classical Bessel function $J_\nu$. We note that the implicit representation of the radii of starlikeness considered in this section were found in \cite{BDM}.

\begin{theorem}\label{th1}
Let $\nu>0.$ Then the radius of starlikeness $r^*\left(f_{\nu}^{(2)}(z;q)\right)$ of the function $$z\mapsto f_{\nu}^{(2)}(z;q)=\left(2^{\nu}c_{\nu}(q)J_{\nu}^{(2)}(z;q)\right)^{\frac{1}{\nu}}$$ is the smallest positive root of the equation $r\cdot dJ_{\nu}^{(2)}(r;q)/dr=0$ and satisfies the following inequality
 $$\frac{4\nu\left(q^{\nu+1}-1\right)\left(q-1\right)}{q^{\nu+1}(\nu+2)}<\left(r^*\left(f_{\nu}^{(2)}(z;q)\right)\right)^2<
 \frac{4\nu(\nu+2)\left(q^{\nu+1}-1\right)\left(q^{\nu+2}-1\right)(q^2-1)}{q^{\nu+1}\left((\nu+2)^2(1+q-q^{\nu+2})-2\nu(\nu+4)q^2+(\nu^2+4\nu-4)q^{\nu+3}\right)}.$$
\end{theorem}

It is worth to mention that by multiplying by $(1-q)^{-2}$ both sides of the above inequality and taking the limit as $q\nearrow1$ for $\nu>0$ we obtain
\begin{equation}\label{eq2.1}
\frac{4\nu(\nu+1)}{\nu+2}<\left(j_{\nu,1}'\right)^2<\frac{4\nu(\nu+1)(\nu+2)^2}{\nu^2+8\nu+8},
\end{equation}
which were obtained earlier by Ismail and Muldoon \cite{ismail3}. Here $j_{\nu,1}'$ is the first positive zero of $J_{\nu}^{\prime}$.

\begin{theorem}\label{th2}
Let $\nu>-1.$ Then the radius of starlikeness $r^*\left(g_{\nu}^{(2)}(z;q)\right)$ of the function $$z\mapsto g_{\nu}^{(2)}(z;q)=2^{\nu}c_{\nu}(q)z^{1-\nu}J_{\nu}^{(2)}(z;q)$$ is the smallest positive root of the equation $r\cdot dJ_{\nu}^{(2)}(r;q)/dr-(\nu-1)J_{\nu}^{(2)}(r;q)=0$ and satisfies the following inequalities
$$2\sqrt{\frac{(q^{\nu+1}-1)(q-1)}{3q^{\nu+1}}}<r^*\left(g_{\nu}^{(2)}(z;q)\right)<
2\sqrt{\frac{3(q^{\nu+1}-1)(q^{\nu+2}-1)(1-q^2)}{q^{\nu+1}(9q-9q^{\nu+2}+q^{\nu+3}-10q^2+9)}}$$
and
$$2\sqrt[4]{\frac{(q^{\nu+1}-1)^2(1-q^{\nu+2})(q-1)^2(q+1)}{q^{2\nu+2}(9q-9q^{\nu+2}+q^{\nu+3}-10q^2+9)}}<
r^*\left(g_{\nu}^{(2)}(z;q)\right)<\frac{2}{3}\sqrt{\frac{3(1-q^{\nu+1})(q^{\nu+3}-1)(q-1)}{(q+1)q^{\nu+1}}{q_1}^*({\nu})},$$
where
$${q_1}^*({\nu})=\frac{(q^3+2q^2+2q+1)(9q-9q^{\nu+2}+q^{\nu+3}-10q^2+9)}{a_{\nu}(q)+b_{\nu}(q)},$$
$$a_{\nu}(q)=-9q^{\nu+2}-12q^{\nu+3}-21q^{\nu+4}+3q^{\nu+5}+6q^{\nu+6}+q^{\nu+7}+9q^{2\nu+5}+3q^{2\nu+6},$$
$$b_{\nu}(q)=3q^{2\nu+7}+q^{2\nu+8}+18q+3q^2-6q^3-15q^4+7q^6+9.$$
\end{theorem}

Note that by multiplying by $(1-q)^{-1}$ both sides of the above inequalities and taking the limit as $q\nearrow1$ for $\nu>-1$ we obtain the first two inequalities of \cite[Theorem 1]{aktas1}, namely:
\begin{equation}\label{eq2.3}
2\sqrt{\frac{\nu+1}{3}}<r^{\star}(\varphi_{\nu})<2\sqrt{\frac{3(\nu+1)(\nu+2)}{4\nu+13}}
\end{equation}
and
\begin{equation}\label{eq2.4}
2\sqrt[4]{\frac{(\nu+1)^2(\nu+2)}{4\nu+13}}<r^{\star}(\varphi_{\nu})<2\sqrt{\frac{(\nu+1)(\nu+3)(4\nu+13)}{2(4\nu^2+26\nu+49)}},
\end{equation}
where $r^{\star}(\varphi_{\nu})$ stands for the radii of starlikeness of the normalized Bessel function $$z\mapsto \varphi_{\nu}(z)=2^{\nu}\Gamma(\nu+1)z^{1-\nu}J_{\nu}(z).$$

\begin{theorem}\label{th3}
Let $\nu>-1.$ Then the radius of starlikeness $r^*\left(h_{\nu}^{(2)}(z;q)\right)$ of the function $$z\mapsto h_{\nu}^{(2)}(z;q)=2^{\nu}c_{\nu}(q)z^{1-\frac{\nu}{2}}J_{\nu}^{(2)}(\sqrt{z};q)$$ is the smallest positive root of the equation $r\cdot dJ_{\nu}^{(2)}(r;q)/dr-(\nu-2)J_{\nu}^{(2)}(r;q)=0$ and satisfies the following inequalities $$\frac{2(1-q)(1-q^{\nu+1})}{q^{\nu+1}}<r^*\left(h_{\nu}^{(2)}(z;q)\right)<
\frac{4(q^{\nu+1}-1)(q^{\nu+2}-1)(1-q^2)}{q^{\nu+1}\left(2q-2q^{\nu+2}+q^{\nu+3}-3q^2+2\right)}$$ and
$$\sqrt{\frac{8(q^{\nu+1}-1)^2(1-q^{\nu+2})(q-1)^2(q+1)}{q^{2\nu+2}\left(2q-2q^{\nu+2}+q^{\nu+3}-3q^2+2\right)}}<
r^*\left(h_{\nu}^{(2)}(z;q)\right)<\frac{4(q^{\nu+1}-1)(q^{\nu+3}-1)(q^3-1)}{q^{\nu+1}}{q_2}^*({\nu}),$$ where $${q_2}^*({\nu})=\frac{2q-2q^{\nu+2}+q^{\nu+3}-3q^2+2}{c_{\nu}(q)+d_{\nu}(q)},$$
$$c_{\nu}(q)=4q^{\nu+2}+3q^{\nu+3}+7q^{\nu+4}-6q^{\nu+5}-5q^{\nu+6}+3q^{\nu+7}-4q^{2\nu+5}+q^{2\nu+6},$$
$$d_{\nu}(q)=q^{2\nu+7}-q^{2\nu+8}-8q+q^2+5q^3+9q^4-6q^6-4.$$
\end{theorem}

Here we would like to emphasize that by multiplying by $(1-q)^{-2}$ both sides of the above inequalities and taking the limit as $q\nearrow1$ for $\nu>-1$ we obtain the first two inequalities of \cite[Theorem 2]{aktas1}, namely:
\begin{equation}\label{eq2.5}
2(\nu+1)<r^{\star}(\phi_{\nu})<\frac{8(\nu+1)(\nu+2)}{\nu+5}
\end{equation}	
and
\begin{equation}\label{eq2.6}
\frac{4(\nu+1)\sqrt{\nu+2}}{\sqrt{\nu+5}}<r^{\star}(\phi_{\nu})<\frac{4(\nu+1)(\nu+3)(\nu+5)}{\nu^2+8\nu+23},
\end{equation}
where $r^{\star}(\phi_{\nu})$ stands for the radii of starlikeness of the normalized Bessel function  $$z\mapsto \phi_{\nu}(z)=2^{\nu}\Gamma(\nu+1)z^{1-\frac{\nu}{2}}J_{\nu}(\sqrt{z}).$$

\begin{theorem}\label{th4}
Let $\nu>0.$ Then the radius of starlikeness $r^*\left(f_{\nu}^{(3)}(z;q)\right)$ of the function $$z\mapsto f_{\nu}^{(3)}(z;q)=\left(c_{\nu}(q)J_{\nu}^{(3)}(z;q)\right)^{\frac{1}{\nu}}$$ is the smallest positive root of the equation $r\cdot dJ_{\nu}^{(3)}(r;q)/dr=0$ and satisfies the following inequality $$\frac{\nu(1-q)(1-q^{\nu+1})}{q(\nu+2)}<\left(r^*\left(f_{\nu}^{(3)}(z;q)\right)\right)^2<
\frac{\nu(\nu+2)(1-q^2)(1-q^{\nu+1})(1-q^{\nu+2})}{q\left((1+q)(1-q^{\nu+2})(\nu+2)^2-2\nu(\nu+4)q(1-q^{\nu+1})\right)}.$$
\end{theorem}

Multiplying by $(1-q)^{-2}$ both sides of the above inequality and taking the limit as $q\nearrow1$ for $\nu>0$ we obtain the following inequality
\begin{equation}\label{eq2.7}
\frac{\nu(\nu+1)}{\nu+2}<\left(j_{\nu,1}^{\prime}\right)^2<\frac{\nu(\nu+1)(\nu+2)^2}{\nu^2+8\nu+8}.
\end{equation}

\begin{theorem}\label{th5}
Let $\nu>-1.$ Then the radius of starlikeness $r^*\left(g_{\nu}^{(3)}(z;q)\right)$ of the function $$z\mapsto g_{\nu}^{(3)}(z;q)=c_{\nu}(q)z^{1-\nu}J_{\nu}^{(3)}(z;q)$$ is the smallest positive root of the equation $r\cdot dJ_{\nu}^{(3)}(r;q)/dr-(\nu-1)J_{\nu}^{(3)}(r;q)=0$ and satisfies the following inequalities
$$\sqrt{\frac{(1-q)(1-q^{\nu+1})}{3q}}<r^*\left(g_{\nu}^{(3)}(z;q)\right)<\sqrt{\frac{3(q^{\nu+1}-1)(q^{\nu+2}-1)(q^2-1)}{q(9q^{\nu+3}-q^{\nu+2}+q-9)}}$$
and
$$\sqrt[4]{\frac{(q^{\nu+1}-1)^2(q-1)^2(q^{\nu+2}-1)(q+1)}{q^2(9q^{\nu+3}-q^{\nu+2}+q-9)}}<r^*\left(g_{\nu}^{(3)}(z;q)\right)<
\sqrt{\frac{(q-1)(q^{\nu+1}-1)(q^{\nu+3}-1)}{3q(q+1)}{q_3}^*({\nu})},$$
where $${q_3}^*({\nu})=\frac{(9q^{\nu+3}-q^{\nu+2}+q-9)(q^3+2q^2+2q+1)}{s_{\nu}(q)+r_{\nu}(q)},$$
$$s_{\nu}(q)=6q^{\nu+2}-12q^{\nu+3}-20q^{\nu+4}-12q^{\nu+5}+6q^{\nu+6}+q^{2\nu+5},$$
$$r_{\nu}(q)=3q^{2\nu+6}+3q^{2\nu+7}+9q^{2\nu+8}+q^3+3q^2+3q+9.$$
\end{theorem}

Multiplying by $(1-q)^{-1}$ both sides of the above inequalities and taking the limit as $q\nearrow1$ for $\nu>-1$ we obtain the following inequalities
\begin{equation}\label{eq2.8}
\sqrt{\frac{\nu+1}{3}}<r^{\star}(\varphi_{\nu})<\sqrt{\frac{3(\nu+1)(\nu+2)}{4\nu+13}}
\end{equation}
and
\begin{equation}\label{eq2.9}
\sqrt[4]{\frac{(\nu+1)^2(\nu+2)}{4\nu+13}}<r^{\star}(\varphi_{\nu})<\sqrt{\frac{(\nu+1)(\nu+3)(4\nu+13)}{2(4\nu^2+26\nu+49)}}.
\end{equation}

\begin{theorem}\label{th6}
Let $\nu>-1.$ Then the radius of starlikeness $r^*\left(h_{\nu}^{(3)}(z;q)\right)$ of the function $$z\mapsto h_{\nu}^{(3)}(z;q)=c_{\nu}(q)z^{1-\frac{\nu}{2}}J_{\nu}^{(3)}(\sqrt{z};q)$$ is the smallest positive root of the equation $r\cdot dJ_{\nu}^{(3)}(r;q)/dr-(\nu-2)J_{\nu}^{(3)}(r;q)=0$ and satisfies the following inequalities $$\frac{(1-q)(1-q^{\nu+1})}{2q}<r^*\left(h_{\nu}^{(3)}(z;q)\right)<\frac{(q^2-1)(q^{\nu+1}-1)(q^{\nu+2}-1)}{q(q-q^{\nu+2}+2q^{\nu+3}-2)}$$
and
$$\frac{(q-1)(q^{\nu+1}-1)}{2q}\sqrt{\frac{2(q^{\nu+2}-1)(q+1)}{q-q^{\nu+2}+2q^{\nu+3}-2}}<r^*\left(h_{\nu}^{(3)}(z;q)\right)<
\frac{(q^{\nu+1}-1)(q^{\nu+3}-1)(1-q)}{q(q+1)}{q_4}^*({\nu}),$$
where $${q_4}^*({\nu})=\frac{(q-q^{\nu+2}+2q^{\nu+3}-2)(q^3+2q^2+2q+1)}
{q+q^2-q^3-4-5q^{\nu+2}+3q^{\nu+3}+10q^{\nu+4}+3q^{\nu+5}-5q^{\nu+6}-q^{2\nu+5}+q^{2\nu+6}+q^{2\nu+7}-4q^{2\nu+8}}.$$
\end{theorem}

Multiplying by $(1-q)^{-2}$ both sides of the above inequalities and taking the limit as $q\nearrow1$ for $\nu>-1$ we obtain the following inequalities
\begin{equation}\label{eq2.10}
\frac{\nu+1}{2}<r^{\star}(\phi_{\nu})<\frac{2(\nu+1)(\nu+2)}{\nu+5}
\end{equation}
and
\begin{equation}\label{eq2.11}
(\nu+1)\sqrt{\frac{\nu+2}{\nu+5}}<r^{\star}(\phi_{\nu})<\frac{(\nu+1)(\nu+3)(\nu+5)}{\nu^2+8\nu+23}.
\end{equation}

It is important to mention that by making a comparison among of above obtained inequalities we have that the left-hand side of \eqref{eq2.7} is weaker than the known result of Ismail and Muldoon, stated as the left-hand side of \eqref{eq2.1}. However, the right-hand side of \eqref{eq2.7} improves the known result of Ismail and Muldoon, stated as the right-hand side of \eqref{eq2.1}. On the other hand, the left-hand side of \eqref{eq2.8} is weaker than the left-hand side of \eqref{eq2.3}, while the right-hand side of \eqref{eq2.8} improves the right-hand side of \eqref{eq2.3} Also, the left-hand side of \eqref{eq2.9} is weaker than the left-hand side of \eqref{eq2.4}, while the right-hand side of \eqref{eq2.9} improves the right-hand side of \eqref{eq2.4}. Finally, the left-hand side of \eqref{eq2.10} is weaker than the left-hand side of \eqref{eq2.5}. However, the right-hand side of \eqref{eq2.10} improves the right-hand side of \eqref{eq2.5}. Also, the left-hand side of \eqref{eq2.11} is weaker than the left-hand side of \eqref{eq2.6}, while the right-hand side of \eqref{eq2.11} improves the right-hand side of \eqref{eq2.6}.

\section{Proofs of main results}
\setcounter{equation}{0}

In this section we are going to present the proofs of our main results.

\begin{proof}[\bf Proof of Theorem \ref{th1}]
In view of  \cite[Theorem 1]{BDM} we know that the radius of starlikeness of the function $f_{\nu}^{(2)}(\cdot;q)$ is the smallest positive root of the equation $r\cdot dJ_{\nu}^{(2)}(r;q)/dr=0$. Also from \cite[eq. (2.2)]{BDM} we have that	
\begin{equation}\label{eq3.1}
\mathcal{B}_{\nu}(z,q)=\frac{2^{\nu}c_{\nu}(q)}{{\nu}z^{\nu-1}}\frac{dJ_{\nu}^{(2)}(z;q)}{dz}=
\prod_{n\geq1}\left(1-\frac{z^2}{\left(j_{\nu,n}^{\prime}(q)\right)^2}\right).
\end{equation}
In addition, by using the definition of the Jackson $q$-Bessel function we obtain
\begin{equation}\label{eq3.2}
\mathcal{B}_{\nu}(z,q)=\frac{2^{\nu}c_{\nu}(q)}{{\nu}z^{\nu-1}}\frac{dJ_{\nu}^{(2)}(z;q)}{dz}=
\sum_{n\geq0}\frac{(-1)^n(2n+\nu)q^{n(n+\nu)}}{{\nu}2^{2n}(q;q)_n(q^{\nu+1};q)_n}z^{2n}.
\end{equation}
Now, taking the logarithmic derivative of \eqref{eq3.1} we get
\begin{equation}\label{eq3.3}
\frac{\mathcal{B}_{\nu}^{\prime}(z,q)}{\mathcal{B}_{\nu}(z,q)}=-2\sum_{k\geq0}\Delta_{k+1}z^{2k+1}, |z|<j_{\nu,1}^{\prime}(q),
\end{equation}
where $\Delta_{k}=\sum_{n\geq1}\left({j_{\nu,n}^{\prime}(q)}\right)^{-2k}$. On the other hand, by considering \eqref{eq3.2} we obtain that
\begin{equation}\label{eq3.4}
\frac{\mathcal{B}_{\nu}^{\prime}(z,q)}{\mathcal{B}_{\nu}(z,q)}=\sum_{n\geq0}\epsilon_nz^{2n+1}\Bigg/\sum_{n\geq0}\mu_nz^{2n},
\end{equation}
where $$\epsilon_n=\frac{(-1)^{n+1}(2n+\nu+2)(n+1)q^{(n+1)(n+\nu+1)}}{\nu2^{2n+1}(q;q)_{n+1}(q^{\nu+1};q)_{n+1}}\ \ \ \mbox{and}\ \ \  \mu_n=\frac{(-1)^{n}(2n+\nu)q^{n(n+\nu)}}{\nu2^{2n}(q;q)_{n}(q^{\nu+1};q)_{n}}.$$
By equating \eqref{eq3.3} and \eqref{eq3.4} and making the Cauchy product we obtain the following Euler -Rayleigh sums $\Delta_{k}=\sum_{n\geq1}\left({j_{\nu,n}^{\prime}(q)}\right)^{-2k}$ in terms of $\nu$ and $q$ for $k\in\left\{1,2\right\}$. Namely,
$$\Delta_{1}=\frac{(\nu+2)q^{\nu+1}}{4\nu(1-q)(1-q^{\nu+1})},$$ and $$\Delta_{2}=\frac{q^{2\nu+2}\left((\nu+2)^2q-2\nu(\nu+4)q^2-(\nu+2)^2q^{\nu+2}+(\nu^2+4\nu-4)q^{\nu+3}+(\nu+2)^2\right)}
{16\nu^2(1-q^{\nu+1})(q^{\nu+2}-1)(1-q)^2(1+q)}.$$
By considering the above Euler-Rayleigh sums and using the Euler-Rayleigh inequalities $$\Delta_{k}^{-\frac{1}{k}}<\left(j_{\nu,1}^{\prime}(q)\right)^2<\frac{\Delta_{k}}{\Delta_{k+1}}$$ for $\nu>-1$ we get the following inequality
$$\frac{4\nu\left(q^{\nu+1}-1\right)\left(q-1\right)}{q^{\nu+1}(\nu+2)}<\left(r^*\left(f_{\nu}^{(2)}(z;q)\right)\right)^2<
\frac{4\nu(\nu+2)\left(q^{\nu+1}-1\right)\left(q^{\nu+2}-1\right)(q^2-1)}{q^{\nu+1}\left((\nu+2)^2(1+q-q^{\nu+2})-2\nu(\nu+4)q^2+(\nu^2+4\nu-4)q^{\nu+3}\right)}.$$
It is possible to have more tighter bounds for the radius of starlikeness of the normalized $q$-Bessel function $f_{\nu}^{(2)}$ for other values of $k$, but it would be quite complicated, this is why we restricted ourselves to the first Euler-Rayleigh inequality.
\end{proof}

\begin{proof}[\bf Proof of Theorem \ref{th2}]
In view of \cite[Theorem 1]{BDM} we know that the radius of starlikeness of the function $g_{\nu}^{(2)}(\cdot;q)$ is the smallest positive root of the equation $r\cdot dJ_{\nu}^{(2)}(r;q)/dr-(\nu-1)J_{\nu}^{(2)}(r;q)=0.$ Now, recall that the zeros $j_{\nu,n}(q), n\in\mathbb{N},$ of the Jackson $q$-Bessel function are all real and simple, according to \cite[Theorem 4.2]{ismail1}. Then, the function $g_{\nu}^{(2)}(\cdot;q)$ belongs to the Laguerre-P\'{o}lya class $\mathcal{LP}$ of real entire functions. Since $\mathcal{LP}$ is closed under differentiation the function $z\mapsto\ {dg_{\nu}^{(2)}(z;q)}/dz$ belongs also to the class $\mathcal{LP}.$ Hence the function $z\mapsto\ {dg_{\nu}^{(2)}(z;q)}/dz$ has only real zeros. Also its growth order $\rho$ is $0$, that is $$\rho=\lim_{n\rightarrow\infty}\frac{n\log{n}}{2n\log2+\log(q;q)_{n}+\log(q^{\nu+1};q)_{n}-\log(2n+1)-n(n+\nu)\log{q}}=0,$$
since as $n\rightarrow\infty$ we have $(q;q)_{n}\rightarrow(q;q)_{\infty}<\infty \text{ and } (q^{\nu+1};q)_{n}\rightarrow(q^{\nu+1};q)_{\infty}<\infty$.
Now, by applying Hadamard's Theorem \cite[p. 26]{Levin} we obtain
$$\phi_{\nu}(z;q)=\frac{dg_{\nu}^{(2)}(z;q)}{dz}=\prod_{n\geq1}\left(1-\frac{z^2}{\alpha_{\nu,n}^2(q)}\right),$$ where ${\alpha_{\nu,n}(q)}$ is the $n$th zero of the function $\phi_\nu(\cdot;q)$. Now, via logarithmic derivation of $\phi_\nu(\cdot;q)$ we obtain
\begin{equation}\label{eq3.5}
\frac{\phi_\nu^{\prime}(z;q)}{\phi_\nu(z;q)}=-2\sum_{k\geq0}\delta_{k+1}z^{2k+1}, |z|<{\alpha_{\nu,1}(q)},
\end{equation}
where $\delta_k=\sum_{n\geq1}\left(\alpha_{\nu,n}(q)\right)^{-2k}$.
Also, by using the infinite sum representation of $\phi_{\nu}$ we get
\begin{equation}\label{eq3.6}
\frac{\phi_\nu^{\prime}(z;q)}{\phi_\nu(z;q)}=\sum_{n\geq0}u_{n}z^{2n+1}\Bigg/\sum_{n\geq0}v_{n}z^{2n},
\end{equation}
where $$u_{n}=\frac{(-1)^{n+1}(2n+2)(2n+3)q^{(n+1)(n+\nu+1)}}{2^{2n+2}(q;q)_{n+1}(q^{\nu+1};q)_{n+1}}\ \ \ \mbox{and}\ \ \  v_{n}=\frac{(-1)^{n}(2n+1)q^{n(n+\nu)}}{2^{2n}(q;q)_{n}(q^{\nu+1};q)_{n}}.$$ By comparing \eqref{eq3.5} and \eqref{eq3.6} and matching all terms with the same degree we have the following Euler-Rayleigh sums $\delta_k=\sum_{n\geq1}\alpha_{\nu,n}^{-2k}(q)$ in terms of $\nu$ and $q$. That is, $$\delta_1=\frac{3q^{\nu+1}}{4(q^{\nu+1}-1)(q-1)},$$
$$\delta_2=-\frac{q^{2\nu+2}(9q-9q^{\nu+2}+q^{\nu+3}-10q^2+9)}{16(q^{\nu+1}-1)^2(q^{\nu+2}-1)(q-1)^2(q+1)}$$ and
$$\delta_3=\frac{3q^{3(\nu+1)}(9q-9q^{\nu+2}+q^{\nu+3}-10q^2+9)}{64(q^{\nu+1}-1)^3(q^{\nu+2}-1)(q^{\nu+3}-1)(q-1)^3q_{1}^*(\nu)}.$$
Now, by considering these Euler-Rayleigh sums in the known Euler-Rayleigh inequalities $$\delta_{k}^{-\frac{1}{k}}<\alpha_{\nu,1}^2(q)<\frac{\delta_{k}}{\delta_{k+1}}$$ for $\nu>-1$ and $k\in\{1,2\}$ we obtain the inequalities of this theorem.
\end{proof}

\begin{proof}[\bf Proof of Theorem \ref{th3}]
By putting $\alpha=0$ in part $\bf c$ of \cite[Theorem 1]{BDM} we have that the radius of starlikeness $r^*\left(h_{\nu}^{(2)}(z;q)\right)$ of the function $h_{\nu}^{(2)}(\cdot;q)$ is the smallest positive root of the equation $r\cdot dJ_{\nu}^{(2)}(r;q)/dr-(\nu-2)J_{\nu}^{(2)}(r;q)=0$.	
It is known that the function
\begin{equation}\label{eq3.7}
\Psi_{\nu}(z;q)=dh_{\nu}^{(2)}(z;q)/dz=\sum_{n\geq0}\frac{(-1)^n(n+1)z^nq^{n(n+\nu)}}{2^{2n}(q;q)_{n}(q^{\nu+1};q)_{n}}
\end{equation}
is an entire function of order $\rho=0$, because
$$\lim_{n\rightarrow\infty}\frac{n\log{n}}{\log(q;q)_{n}+\log(q^{\nu+1};q)_{n}+2n\log2-n(n+\nu)\log{q}-\log(n+1)}=0$$ since as $n\rightarrow\infty$ we have $(q;q)_{n}\rightarrow(q;q)_{\infty}<\infty \text{ and } (q^{\nu+1};q)_{n}\rightarrow(q^{\nu+1};q)_{\infty}<\infty$. Moreover, we know that the zeros $\beta_{\nu,n}(q), n\in\mathbb{N},$ of the function $z\mapsto{dh_{\nu}^{(2)}(z;q)}/{dz}$ are real when $\nu>-1$, according to \cite[Lemma 6]{BDM}. Now, by applying Hadamard's Theorem \cite[p. 26]{Levin} we get
\begin{equation}\label{eq3.8}
\Psi_{\nu}(z;q)=dh_{\nu}^{(2)}(z;q)/dz=\prod_{n\geq1}\left(1-\frac{z}{\beta_{\nu,n}^2(q)}\right).
\end{equation}
From \eqref{eq3.7} we have
\begin{equation}\label{eq3.9}
\frac{\Psi_{\nu}^{\prime}(z;q)}{\Psi_{\nu}(z;q)}=\sum_{n\geq0}t_{n}z^n\Bigg/\sum_{n\geq0}m_{n}z^n,
\end{equation}
where $$t_{n}=\frac{(-1)^{n+1}(n+2)(n+1)q^{(n+1)(n+\nu+1)}}{2^{2n+2}(q;q)_{n+1}(q^{\nu+1};q)_{n+1}}\ \ \ \mbox{and}\ \ \  m_{n}=\frac{(-1)^{n}(n+1)q^{n(n+\nu}}{2^{2n}(q;q)_{n}(q^{\nu+1};q)_{n}}.$$
Also, taking the derivative of \eqref{eq3.8} logarithmically we get
\begin{equation}\label{eq3.10}
\frac{\Psi_{\nu}^{\prime}(z;q)}{\Psi_{\nu}(z;q)}=-\sum_{k\geq0}\Theta_{k+1}z^k, \rvert{z}\rvert<\beta_{\nu,1}^2(q),
\end{equation}
where $\Theta_{k}=\sum_{n\geq1}\beta_{\nu,n}^{-2k}(q)$. Now, it is possible to state the Euler-Rayleigh sums $\Theta_{k}=\sum_{n\geq1}\beta_{\nu,n}^{-2k}(q)$ in terms of $\nu$ and $q$. By comparison of the coefficients of \eqref{eq3.9} and \eqref{eq3.10} we obtain
$$\Theta_{1}=\frac{q^{\nu+1}}{2(1-q)(1-q^{\nu+1})}\ \ \ \mbox{and}\ \ \ \Theta_{2}=\frac{q^{2\nu+2}(2q-2q^{\nu+2}+q^{\nu+3}-3q^2+2)}{8(q^{\nu+1}-1)(1-q^{\nu+2})(q-1)^2(q+1)}.$$
By considering the Euler-Rayleigh inequalities $$\Theta_{k}^{-\frac{1}{k}}<r^*\left(h_{\nu}^{(2)}(z;q)\right)<\frac{\Theta_{k}}{\Theta_{k+1}}$$ for $\nu>-1$ and $k\in\{1,2\}$ we get the following inequalities
$$\frac{2(1-q)(1-q^{\nu+1})}{q^{\nu+1}}<r^*\left(h_{\nu}^{(2)}(z;q)\right)<
\frac{4(q^{\nu+1}-1)(q^{\nu+2}-1)(1-q^2)}{q^{\nu+1}\left(2q-2q^{\nu+2}+q^{\nu+3}-3q^2+2\right)}$$ and
$$\sqrt{\frac{8(q^{\nu+1}-1)^2(1-q^{\nu+2})(q-1)^2(q+1)}{q^{2\nu+2}\left(2q-2q^{\nu+2}+q^{\nu+3}-3q^2+2\right)}}
<r^*\left(h_{\nu}^{(2)}(z;q)\right)<\frac{4(q^{\nu+1}-1)(q^{\nu+3}-1)(q^3-1)}{q^{\nu+1}}{q_2}^*({\nu}).$$
\end{proof}

\begin{proof}[\bf Proof of Theorem \ref{th4}]
By taking $\alpha=0$ in part $\bf a$ of \cite[Theorem 1]{BDM} we have that the radius of starlikeness $r^*\left(f_{\nu}^{(3)}(z;q)\right)$ of the function $f_{\nu}^{(3)}(\cdot;q)$ is the smallest positive root of the equation $r\cdot dJ_{\nu}^{(3)}(r;q)/dr=0$. Also, from the equation $(2.2)$ in \cite{BDM} we have that	
\begin{equation}\label{eq3.11}
\mathcal{K}_{\nu}(z,q)=\frac{c_{\nu}(q)}{{\nu}z^{\nu-1}}\frac{dJ_{\nu}^{(3)}(z;q)}{dz}=
\prod_{n\geq1}\left(1-\frac{z^2}{\left(l_{\nu,n}^{\prime}(q)\right)^2}\right).
\end{equation}
On the other hand, by using the infinite sum representation of the Hahn-Exton $q$-Bessel function we obtain
\begin{equation}\label{eq3.12}
\mathcal{K}_{\nu}(z,q)=\frac{c_{\nu}(q)}{{\nu}z^{\nu-1}}\frac{dJ_{\nu}^{(3)}(z;q)}{dz}=
\sum_{n\geq0}\frac{(-1)^n(2n+\nu)q^{\frac{1}{2}n(n+1)}}{{\nu}(q;q)_n(q^{\nu+1};q)_n}z^{2n}.
\end{equation}	
Taking the logarithmic derivative of \eqref{eq3.11} we have
\begin{equation}\label{eq3.13}
\frac{\mathcal{K}_{\nu}^{\prime}(z,q)}{\mathcal{K}_{\nu}(z,q)}=-2\sum_{k\geq0}\xi_{k+1}z^{2k+1}, |z|<{l_{\nu,1}^{\prime}}(q),
\end{equation}
where $\xi_k=\sum_{n\geq}\left(l_{\nu,n}^{\prime}(q)\right)^{-2k}$. By using the expression \eqref{eq3.12} we write
\begin{equation}\label{eq3.14}
\frac{\mathcal{K}_{\nu}^{\prime}(z,q)}{\mathcal{K}_{\nu}(z,q)}=\sum_{n\geq0}s_{n}z^{2n+1}\Bigg/\sum_{n\geq0}p_{n}z^{2n},
\end{equation}
where $$s_{n}=\frac{(-1)^{n+1}2(n+1)(2n+\nu+2)q^{\frac{1}{2}(n+1)(n+2)}}{\nu(q;q)_{n+1}(q^{\nu+1};q)_{n+1}}\ \ \ \mbox{and}\ \ \  p_n=\frac{(-1)^{n}(2n+\nu)q^{\frac{1}{2}n(n+1)}}{\nu(q;q)_{n}(q^{\nu+1};q)_{n}}.$$ Now, it is possible to express the Euler-Rayleigh sums $\xi_k=\sum_{n\geq}\left(l_{\nu,n}^{\prime}(q)\right)^{-2k}$ in terms of $\nu$ and $q$. By matching the coefficients of the equalities \eqref{eq3.13} and \eqref{eq3.14} we get
$$\xi_{1}=\frac{q(\nu+2)}{\nu(1-q)(1-q^{\nu+1})}\ \ \ \mbox{and}\ \ \ \xi_{2}=q^2\left(\frac{(1+q)(1-q^{\nu+2})(\nu+2)^2-2\nu(\nu+4)q(1-q^{\nu+1})}{\nu^2(1-q)^2(1+q)(1-q^{\nu+1})^2(1-q^{\nu+2})}\right).$$ By considering the above Euler-Rayleigh sums in the Euler-Rayleigh inequalities $$\xi_{k}^{-\frac{1}{k}}<r^*\left(f_{\nu}^{(3)}(z;q)\right)<\frac{\xi_{k}}{\xi_{k+1}}$$ for $\nu>0$ and $k=1$ we have the following inequality
$$\frac{\nu(1-q)(1-q^{\nu+1})}{q(\nu+2)}<r^*\left(f_{\nu}^{(3)}(z;q)\right)^2<
\frac{\nu(\nu+2)(1-q^2)(1-q^{\nu+1})(1-q^{\nu+2})}{q\left((1+q)(1-q^{\nu+2})(\nu+2)^2-2\nu(\nu+4)q(1-q^{\nu+1})\right)}.$$
\end{proof}

\begin{proof}[\bf Proof of Theorem \ref{th5}]
By virtue of part {\bf b} of \cite[Theorem 1]{BDM} the radius of starlikeness $r^*\left(g_{\nu}^{(3)}(z;q)\right)$ of the function $g_{\nu}^{(3)}(\cdot;q)$ is the smallest positive root of the equation $r\cdot dJ_{\nu}^{(3)}(r;q)/dr-(\nu-1)J_{\nu}^{(3)}(r;q)=0$.
Now, recall that the zeros $j_{\nu,n}(q), n\in\mathbb{N},$ of the Hahn-Exton $q$-Bessel function are all real and simple, according to \cite[Theorem 4.2]{ismail1}. Then, the function $g_{\nu}^{(3)}(\cdot;q)$ belongs to the Laguerre-P\'{o}lya class $\mathcal{LP}$ of real entire functions. Since $\mathcal{LP}$ is closed under differentiation the function
\begin{equation}\label{eq3.15}
z\mapsto\Phi_{\nu}(z;q)=\frac{dg_{\nu}^{(3)}(z;q)}{dz}=\sum_{n\geq0}\frac{(-1)^n(2n+1)q^{\frac{1}{2}n(n+1)}}{(q;q)_{n}(q^{\nu+1};q)_{n}}z^{2n}
\end{equation}
belongs also to the class $\mathcal{LP}$. Hence the function $z\mapsto{dg_{\nu}^{(3)}(z;q)}/{dz}$ has only real zeros. Also its growth order $\rho$ is $0$, that is $$\rho=\lim_{n\rightarrow\infty}\frac{n\log{n}}{\log(q;q)_{n}+\log(q^{\nu+1};q)_{n}-\frac{1}{2}n(n+1)\log{q}-\log(2n+1)}=0,$$
since as $n\rightarrow\infty$ we have $(q;q)_{n}\rightarrow(q;q)_{\infty}<\infty \text{ and } (q^{\nu+1};q)_{n}\rightarrow(q^{\nu+1};q)_{\infty}<\infty$. Due to Hadamard's Theorem \cite[p. 26]{Levin} we have
\begin{equation}\label{eq3.16}
\Phi_{\nu}(z;q)=\prod_{n\geq}\left(1-\frac{z^2}{\gamma_{\nu,n}^2(q)}\right).
\end{equation}	
If we consider the equality \eqref{eq3.15} we can write that
\begin{equation}\label{eq3.17}
\frac{\Phi_{\nu}^{\prime}(z;q)}{\Phi_{\nu}(z;q)}=\sum_{n\geq0}k_{n}z^{2n+1}\Bigg/\sum_{n\geq0}l_{n}z^{2n},
\end{equation}
where $$k_{n}=\frac{(-1)^{n+1}(2n+2)(2n+3)q^{\frac{1}{2}(n+1)(n+2)}}{(q;q)_{n+1}(q^{\nu+1};q)_{n+1}}\ \ \ \mbox{and} \ \ \ l_{n}=\frac{(-1)^{n}(2n+1)q^{\frac{1}{2}n(n+1)}}{(q;q)_{n}(q^{\nu+1};q)_{n}}.$$ In addition, by using the logarithmic derivative of both sides of \eqref{eq3.16} we have
\begin{equation}\label{eq3.18}
\frac{\Phi_{\nu}^{\prime}(z;q)}{\Phi_{\nu}(z;q)}=-2\sum_{k\geq0}S_{k+1}z^{2k+1}, \lvert{z}\lvert<\gamma_{\nu,1}(q),
\end{equation}
where $S_{k}=\sum_{n\geq1}\gamma_{\nu,n}^{-2k}(q)$.
Now, it is possible to express the Euler-Rayleigh sums $S_{k}=\sum_{n\geq1}\gamma_{\nu,n}^{-2k}(q)$ in terms of $\nu$ and $q$. By matching the coefficients of the equalities \eqref{eq3.17} and \eqref{eq3.18} we get
$$S_{1}=\frac{3q}{(1-q)(1-q^{\nu+1})}, S_{2}=\frac{q^2(9q^{\nu+3}-q^{\nu+2}+q-9)}{(q^{\nu+1}-1)^2(q^{\nu+2}-1)(q-1)^2(q+1)}$$ and $$S_{3}=\frac{3q^3(9q^{\nu+3}-q^{\nu+2}+q-9)}{(q^{\nu+1}-1)^3(q^{\nu+2}-1)(q^{\nu+3}-1)(q-1)^3{q_3}^*({\nu})}.$$
By considering the above Euler-Rayleigh sums in the Euler-Rayleigh inequalities $$S_{k}^{-\frac{1}{k}}<\left(r^*\left(g_{\nu}^{(3)}(z;q)\right)\right)^2<\frac{S_{k}}{S_{k+1}}$$ for $\nu>-1$ and $k\in\{1,2\}$ we have the following inequalities
$$\sqrt{\frac{(1-q)(1-q^{\nu+1})}{3q}}<r^*\left(g_{\nu}^{(3)}(z;q)\right)<\sqrt{\frac{3(q^{\nu+1}-1)(q^{\nu+2}-1)(q^2-1)}{q(9q^{\nu+3}-q^{\nu+2}+q-9)}}$$
and
$$\sqrt[4]{\frac{(q^{\nu+1}-1)^2(q-1)^2(q^{\nu+2}-1)(q+1)}{q^2(9q^{\nu+3}-q^{\nu+2}+q-9)}}<r^*\left(g_{\nu}^{(3)}(z;q)\right)<\sqrt{\frac{(q-1)(q^{\nu+1}-1)(q^{\nu+3}-1)}{3q(q+1)}{q_3}^*({\nu})}.$$
\end{proof}

\begin{proof}[\bf Proof of Theorem \ref{th6}]
Thanks to part {\bf c} of \cite[Theorem 1]{BDM} we know that the radius of starlikeness $r^*\left(h_{\nu}^{(3)}(z;q)\right)$ of the function $h_{\nu}^{(3)}(\cdot;q)$ is the smallest positive root of the equation $r\cdot dJ_{\nu}^{(3)}(r;q)/dr-(\nu-2)J_{\nu}^{(3)}(r;q)=0$.
On the other hand, with the help of the infinite sum representation of the Hahn-Exton $q$-Bessel function the function $z\mapsto{dh_{\nu}^{(3)}(z;q)}/{dz}$ can be written as an infinite sum as follow:
\begin{equation}\label{eq3.19}
\Upsilon_{\nu}(z;q)=\frac{dh_{\nu}^{(3)}(z;q)}{dz}=\sum_{n\geq0}\frac{(-1)^n(n+1)q^{\frac{1}{2}n(n+1)}}{(q;q)_{n}(q^{\nu+1};q)_{n}}z^n.
\end{equation}
Also its growth order $\rho$ is $0$, that is $$\rho=\lim_{n\rightarrow\infty}\frac{n\log{n}}{\log(q;q)_{n}+\log(q^{\nu+1};q)_{n}-\frac{1}{2}n(n+1)\log{q}-\log(n+1)}=0,$$
since as $n\rightarrow\infty$ we have $(q;q)_{n}\rightarrow(q;q)_{\infty}<\infty \text{ and } (q^{\nu+1};q)_{n}\rightarrow(q^{\nu+1};q)_{\infty}<\infty$. In addition, it is known that the zeros $\gamma_{\nu,n}(q)$ of the function $\Upsilon_{\nu}(z;q)$ are real for $\nu>-1$ and $n\in\mathbb{N}$, according to \cite[Lemma 6]{BDM}. By applying Hadamard's Theorem \cite[p. 26]{Levin} we have
\begin{equation}\label{eq3.20}
\Upsilon_{\nu}(z;q)=\prod_{n\geq1}\left(1-\frac{z}{\gamma_{\nu,n}^{2}(q)}\right).
\end{equation}
Now, using the equality \eqref{eq3.19} we have
\begin{equation}\label{eq3.21}
\frac{\Upsilon_{\nu}^{\prime}(z;q)}{\Upsilon_{\nu}(z;q)}=\sum_{n\geq0}x_{n}z^n\Bigg/\sum_{n\geq0}y_{n}z^n,
\end{equation}
where $$x_{n}=\frac{(-1)^{n+1}(n+1)(n+2)q^{\frac{1}{2}(n+1)(n+2)}}{(q;q)_{n+1}(q^{\nu+1};q)_{n+1}}\text{ and }y_{n}=\frac{(-1)^{n}(n+1)q^{\frac{1}{2}n(n+1)}}{(q;q)_{n}(q^{\nu+1};q)_{n}}.$$ Then, taking the logarithmic derivative of \eqref{eq3.20} we get
\begin{equation}\label{eq3.22}
\frac{\Upsilon_{\nu}^{\prime}(z;q)}{\Upsilon_{\nu}(z;q)}=-\sum_{k\geq0}\sigma_{k+1}z^k, \rvert{z}\rvert<\gamma_{\nu,1}^2(q),
\end{equation}
where $\sigma_{k}=\sum_{n\geq1}\gamma_{\nu,n}^{-2k}(q)$. By comparing all coefficients of \eqref{eq3.21} and \eqref{eq3.22} it is possible to express the Euler-Rayleigh sums $\sigma_{k}=\sum_{n\geq1}\gamma_{\nu,n}^{-2k}(q)$ in terms of $\nu\text{ and } q.$ Namely, $$\sigma_{1}=\frac{6q^3}{(q-1)(1-q^2)(1-q^{\nu+1})(1-q^{\nu+2})}, \sigma_{2}=\frac{2q^2(q-q^{\nu+2}+2q^{\nu+3}-2)}{(q^{\nu+1}-1)^2(q^{\nu+2}-1)(q-1)^2(q+1)}$$ and $$\sigma_{3}=\frac{2q^3(q-q^{\nu+2}+2q^{\nu+3}-2)}{(q^{\nu+1}-1)^3(q^{\nu+2}-1)(1-q^{\nu+3})(q-1)^3q_{4}^{*}(\nu)}.$$ Now, applying the Euler-Rayleigh inequalities $$\sigma_{k}^{-\frac{1}{k}}<r^*\left(h_{\nu}^{(3)}(z;q)\right)<\frac{\sigma_{k}}{\sigma_{k+1}}$$ for $k\in\{1,2\}$ we get the following inequalities
$$\frac{(q-1)(1-q^2)(1-q^{\nu+1})(1-q^{\nu+2})}{6q^3}<r^*\left(h_{\nu}^{(3)}(z;q)\right)<
\frac{(q^2-1)(q^{\nu+1}-1)(q^{\nu+2}-1)}{q(q-q^{\nu+2}+2q^{\nu+3}-2)}$$ and
$$\frac{(q-1)(q^{\nu+1}-1)}{2q}\sqrt{\frac{2(q^{\nu+2}-1)(q+1)}{q-q^{\nu+2}+2q^{\nu+3}-2}}<
r^*\left(h_{\nu}^{(3)}(z;q)\right)<\frac{(q^{\nu+1}-1)(q^{\nu+3}-1)(1-q)}{q(q+1)}{q_4}^*({\nu}).$$
\end{proof}

\end{document}